\documentclass[12pt]{article}

\usepackage[top=2.4cm,bottom=2.2cm,left=2.6cm,right=2cm]{geometry}
\usepackage{latexsym,epsfig,psfrag}
\usepackage{eepic,colordvi,amscd}

\usepackage{amsmath,amsthm,dsfont,amsfonts,amssymb,fancyhdr}
\usepackage{enumerate}
\usepackage[usenames,dvipsnames]{color}
\usepackage{bbm,soul,supertabular,longtable,verbatim,extarrows,xcolor}
\usepackage{titlesec}
\usepackage{BOONDOX-cal}
\usepackage{cite}
\usepackage{tikz}
\usepackage{booktabs,multirow,makecell}
\usepackage{authblk}
\usepackage{color}
\usetikzlibrary[trees]
\usepackage {url}

\newcommand{\SPLS}{SPLS$(\sigma)$ }

\newcommand{\bal}{\ensuremath{t\alpha+\alpha-1}}

\usepackage{lipsum}
\usepackage{mathrsfs}
\usepackage{indentfirst}
\usepackage[colorlinks=true, allcolors=blue]{hyperref}
\newtheorem{thm}{Theorem}[section]
\newtheorem{cor}[thm]{Corollary}
\newtheorem{lem}[thm]{Lemma}

\theoremstyle{definition}

\theoremstyle{remark}
\newtheorem{rem}[thm]{Remark}
\numberwithin{equation}{subsection}
\newcounter{si}
\addtocounter{si}{1}

\usepackage{mathtools}

\newcommand{\Rmnum}[1]{\uppercase\expandafter{\romannumeral #1}}  

\makeatletter
\newcommand\blfootnote[1]{%
  \begingroup
  \renewcommand\thefootnote{}\footnote{#1}%
  \addtocounter{footnote}{-1}%
  \endgroup
}
\makeatother

\begin{document}
\title{An improved bound for strongly regular graphs with smallest eigenvalue $-m$}

\author[a,b]{Jack H. Koolen}
\author[a]{Chenhui Lv}
\author[c]{Greg Markowsky}
    \author[d]{Jongyook Park\footnote{J. Park is the corresponding author.}}

 \affil[a]{\footnotesize{School of Mathematical Sciences, University of Science and
Technology of China, Hefei, 230026, People's Republic of China}}
\affil[b]{\footnotesize{CAS Wu Wen-Tsun Key Laboratory of Mathematics, University
of Science and Technology of China, Hefei, 230026, People's Republic of China}}
\affil[c]{\footnotesize{Department of Mathematics,  Monash University, Australia}} 
\affil[d]{\footnotesize{Department of Mathematics, Kyungpook National University, Daegu, 41566, Republic of Korea}}


\maketitle

\blfootnote{\small E-mail addresses: {\tt koolen@ustc.edu.cn} (J.H. Koolen),  {\tt lch1994@mail.ustc.edu.cn} (C. Lv),  {\tt gmarkowsky@gmail.com} (G. Markowsky),  {\tt jongyook@knu.ac.kr} (J. Park)}
\begin{abstract}
In 1979, Neumaier gave a bound on $\lambda$ in terms of $m$ and $\mu$, where $-m$ is the smallest eigenvalue of a primitive strongly regular graph, unless the graph in question belongs to one of the two infinite families of strongly regular graphs. We improve this result. We also indicate how our methods can be used to give an alternate derivation of Bruck's Completion Theorem for orthogonal arrays.
\end{abstract}


{\bf AMS 2020 MSC:} 05B15, 05C50, 05C60, 05E30

\section{Introduction} \label{intro}

In this paper, we study strongly regular graphs with a fixed smallest eigenvalue (for the definition, see Section~\ref{prelim}). Sims (see~\cite[Theorem~8.6.4]{SRG}) proved the following:

\begin{thm}[cf.~{\cite[Theorem 8.6.4]{SRG}}]\label{srgsims}
Let $\Gamma$ be a primitive strongly regular graph with smallest eigenvalue $-m$, where $m$ is a positive integer. Then there exists a finite list $\mathcal{L}(m)$ such that one of the following holds:

\begin{enumerate}[(i)]
    \item[(i)] $\Gamma$ is a Latin square graph.
    \item[(ii)] $\Gamma$ is the block graph of a Steiner system.
    \item[(iii)] $\Gamma$ is a member of $\mathcal{L}(m)$.
\end{enumerate}
\end{thm}

Neumaier, building on results of Bruck and Bose, made the statement of Sims explicit by completing the following theorem.

\begin{thm}[cf.~{\cite[Theorem 8.6.3]{SRG}}]\label{srgneumaier}
Let $\Gamma$ be a primitive strongly regular graph with smallest eigenvalue $-m$, where $m$ is a positive integer. Let
$f(m, \mu) = \frac{1}{2} m (m - 1) (\mu + 1) + \mu - m - 1$.Then the following hold:
\begin{enumerate}
\item[(i)] (Bruck \cite{bruck1963finite}) If $\mu = m(m - 1)$ and $\lambda > f(m, \mu)$, then $\Gamma$ is the
collinearity graph of a partial geometry $\mathrm{pg}(s, t, \alpha)$ with $\alpha = t$; that is, a
Latin square graph $\mathrm{LS}_m(\lambda - m(m - 3))$.
\item[(ii)] (Bose \cite{bose1963}) If $\mu = m^2$ and $\lambda > f(m, \mu)$, then $\Gamma$ is the
collinearity graph of a partial geometry $\mathrm{pg}(s, t, \alpha)$ with $\alpha = t+1$; that is, the block
graph of a $2$-$(\lambda(m - 1) - m(m - 1)(m - 2) + m, m, 1)$ design.
\item[(iii)] (Neumaier \cite[Theorem 4.7]{neumaier-m}) If $\mu \neq m(m - 1)$ and $\mu \neq m^2$, then $\lambda \leq f(m,
\mu)$.
\end{enumerate}
\end{thm}

Note that Neumaier also showed that for a primitive strongly regular graph with smallest eigenvalue $-m$, the parameter $\mu$ is bounded by $(2m-3)m^3$(see {\cite[Theorem 3.1]{neumaier-m}}). These two facts imply that the second largest eigenvalue is
bounded by a function in $m$, and hence the valency is also bounded by a function in $m$, unless the graph belongs
to one of the two infinite families mentioned in the theorem.

In this paper, we improve the result of Neumaier as follows.

\begin{thm}\label{mainclassification}
Let $\Gamma$ be a primitive strongly regular graph with smallest
eigenvalue $-m$, where $m$ is a positive integer. Let $f(m, \mu) = \frac{8}{3} m (\mu
- 1) - \frac{2}{3} \mu +3m  - \frac{10}{3}$. Then the following hold:
\begin{enumerate}
\item[(i)] If $\mu = m(m - 1)$ and $\lambda > f(m, \mu)$, then
$\Gamma$ is the collinearity graph of a partial geometry
$\mathrm{pg}(s, t, \alpha)$ with $\alpha = t$; that is, a Latin square
graph $\mathrm{LS}_m(\lambda - m(m - 3))$.
\item[(ii)]If $\mu = m^2$ and $\lambda > f(m, \mu)$, then
$\Gamma$ is the collinearity graph of a partial geometry
$\mathrm{pg}(s, t, \alpha)$ with $\alpha = t+1$; that is, the block graph of a
$2$-$(\lambda(m - 1) - m(m - 1)(m - 2) + m, m, 1)$ design.
\item[(iii)] If $\mu \neq m(m - 1)$ and $\mu \neq m^2$, then
$\lambda \leq f(m, \mu)$.
\end{enumerate}
\end{thm}

We should mention that Metsch \cite{Metsch} has already shown Theorem \ref{mainclassification} under the additional condition that the strongly regular graph in question has geometric parameters lying within a particular range of values. Without this additional condition, Metsch's method and result are not valid. Our contribution is showing that this additional condition is not required, and giving a more elementary and general method for proving the result. We have endeavored to keep the paper as simple as possible; as a consequence, our result is slightly weaker than Metsch's when the graph in question does satisfy the additional requirement. However, with a bit of extra effort (which we do not expend in this paper) our methods can be adapted to prove exactly Metsch's result in the case that he proved it.

In comparison to Neumaier's result (Theorem \ref{srgneumaier}), our main result (Theorem \ref{mainclassification}) provides a better bound when $m \geq 6$. For $m \leq 5$, Neumaier's bound is typically better, except in a few exceptional cases.

We now give an outline of the proof. The condition $\lambda>f(m,\mu)$ ensure that the strongly regular graph is the point graph of a partial linear space, using a result of Metsch (see Theorem \ref{metsch}).
Using the line-vertex incidence matrix, we show that a strongly regular graph satisfying the condition $\lambda>f(m,\mu)$ has parameters which would match those of the point graph of a partial geometry. Then we use the method of \cite{guo2020nonexistence} to show that our strongly regular graph must in fact be geometric (that is, it must be the point graph of a partial geometry). Finally, Theorem \ref{mainclassification} follows easily from this.

The following are several applications of Neumaier's result (Theorem~\ref{srgneumaier}) and our main result (Theorem~\ref{mainclassification}): (1) In the appendix, we demonstrate how Theorem~\ref{mainclassification}~(i) can be applied to provide a proof of Bruck's well-known Completion Theorem;
(2) The class of strongly regular graphs is known to be challenging for the graph isomorphism problem; see, for example, Babai~\cite{Babai2014,Babai2015}.
Another application of Neumaier's result was given by Spielman~\cite{Spielman1996} and later improved by Babai et~al.~\cite{Babai2013Faster}, who used it to derive bounds on the complexity of the graph isomorphism problem for the class of strongly regular graphs. For both results, providing a good bound on $\lambda$ is crucial.
Spielman's bound on $\lambda$~\cite{Spielman1996} was subsequently improved by Pyber~\cite{Pyber2014} and later by Babai and Wilmes~\cite{BabaiWilmes2016}. 
Our bound, stated in Theorem~\ref{mainclassification}, improves upon that of Babai and Wilmes (see Section~\ref{compBWBound}) and is therefore likely to contribute to a deeper understanding of this important question; 
(3) in \cite{Pyber2014}, Pyber showed that, apart from finitely many exceptions, all connected strongly regular graphs are Hamiltonian, and that Neumaier's result plays a crucial role in his proof; 
(4) Neumaier's result has further applications. For example, in \cite{koolen1homo}, the authors classify $1$-homogeneous distance-regular graphs of large valency using Theorem~\ref{srgneumaier}.

The paper is organized as follows. The next section reviews (without proofs) the parts of the theory of strongly regular graphs, partial linear spaces, and partial geometries that are relevant to our study. We introduce generic parameters, geometric parameters, and Delsarte cliques, followed by a brief discussion of the relationships between them.  Section \ref{sec: partial lin space} applies many of the methods developed in \cite{guo2020nonexistence} in a much more general setting, and it is revealed there that these methods really are about the presence or absence of Delsarte cliques in a graph; this interpretation was not noted in \cite{guo2020nonexistence}. In Section \ref{sec: proof}, we finally are in a position to prove our main result. In Section~\ref{compBWBound}, we compare our bound with that of Babai and Wilmes. We have also included an appendix, Appendix \ref{sec: orth array}, which deals with the connection between our results and Bruck's well-known completion theorem for orthogonal arrays.

\section{Preliminaries} \label{prelim}

All graphs under consideration are assumed to be simple and connected, unless otherwise stated. We also assume graphs are not complete, to avoid the constant need to qualify statements. The following are standard notations. Let $\Gamma$ be a graph with vertex set $V(\Gamma)$. For two vertices $x, y\in V(\Gamma)$, $d(x,y)$ is the {\it distance} between them. The {\it diameter} $D$ of $\Gamma$ is the maximum distance among two vertices of $\Gamma$. For a vertex $x$, we let $\Gamma_i(x) = \{y \in V(\Gamma): d(x,y)=i\}$ $(0\leq i\leq D)$, and write $\Gamma(x) := \Gamma_1(x)$. The {\it eigenvalues} of $\Gamma$ are the eigenvalues of its adjacency matrix.

\subsection{Strongly regular graphs}
An {\it amply regular graph} with parameters $(v, k, \lambda, \mu)$ is a $k$-regular graph on
$v$ vertices, with the property that every pair of adjacent vertices have exactly $\lambda$
common neighbors and every pair of vertices at distance two have exactly $\mu$ common
neighbors. A {\it strongly regular graph} with {\it standard parameters} $(v, k, \lambda, \mu)$ is an amply regular graph with the same parameters and diameter $2$, so that each pair of distinct vertices have either $\lambda$ or $\mu$ common neighbors, depending on whether or not they are adjacent.  Additional information on strongly regular graphs and their properties can be found in \cite{SRG}.  To dispense with trivialities, we will be concerned only with {\it primitive} strongly regular graphs, which are connected strongly regular graphs whose complement is also connected. 

An important property of strongly regular graphs is that they have exactly $3$ distinct eigenvalues, among which the valency $k$ is always one of them. The eigenvalues different from the valency are called {\it nontrivial}. An essential fact about strongly regular graphs is that the nontrivial eigenvalues $\theta_1$ and $\theta_2$ are the positive and negative roots of the quadratic equation $x^2 - (\lambda - \mu)x - (k - \mu)=0$. Consequently, we have the following fundamental relations:
\begin{equation} \label{eigeq}
-\theta_1 \theta_2 = k - \mu \quad \text{and} \quad \theta_1 + \theta_2 = \lambda - \mu.
\end{equation}

By a simple counting argument, we obtain $k(k-\lambda-1) = (v-k-1)\mu$, and hence $v$ is determined by $k$, $\lambda$, and $\mu$.

Strongly regular graphs with standard parameters $(v, k, \lambda, \mu) = (4t + 1, 2t, t - 1, t)$ are
known as \emph{conference graphs}. It is well known that  all eigenvalues of a strongly regular graph are integers, unless the graph is a conference graph (see, for example,  \cite[Lemma~10.3.3]{GTM207}). Throughout this paper, we always assume that $\Gamma$ is not a conference graph.

We now introduce the notion of  generic parameters, which extends the concept of geometric parameters to the real domain. The geometric parameters themselves will be discussed in Section~\ref{sec:partialgeo}. This extension is advantageous because every strongly regular graph admits (unique) generic parameters.

Let $\Gamma$ be a strongly regular graph with standard parameters $(v, k, \lambda, \mu)$. We say that $\Gamma$ has \emph{generic parameters} $(s, t, \alpha)$ if there exist real numbers $s, t, \alpha > 0$ that satisfy the following conditions:
\begin{equation} \label{standparsrg}
\begin{gathered} 
\hspace{-0.2cm}k = s(t+1),\\
\hspace{1.2cm}\lambda = s + t(\alpha - 1 ) - 1, \\
\hspace{-0.25cm}\mu = (t + 1)\alpha. \\
\end{gathered}    
\end{equation}
Furthermore, by equation~\eqref{eigeq}, the nontrivial eigenvalues of $\Gamma$ are $s - \alpha$ and $-t - 1$.

On the other hand, let $\theta_{\min}$ be the smallest eigenvalue of $\Gamma$. We may define
\begin{equation} \label{classparsrg}
\begin{gathered} 
\hspace{-0.6cm}s = \frac{k}{-\theta_{\min}},\\ 
\hspace{0.1cm}t = -\theta_{\min} - 1, \\ 
\hspace{-0.6cm}\alpha = \frac{\mu}{-\theta_{\min}}, 
\end{gathered}    
\end{equation}
so that, again by equation~\eqref{eigeq}, $\Gamma$ has generic parameters $(s, t, \alpha)$. Therefore, each strongly regular graph admits a unique set of generic parameters $(s,t,\alpha)$, and one may freely pass between the standard parameters and the generic parameters. Moreover, since we study strongly regular graphs with fixed smallest eigenvalue, it is advantageous to include the smallest eigenvalue in the parameter set.

\subsection{Delsarte cliques} \label{sec: classpar}
In certain classes of $k$-regular graphs, the bound $|C| \leq 1 + \frac{k}{-\theta_{\min}}$ must be obeyed, where $C$ is any clique and $\theta_{\min}$ is the smallest eigenvalue of the graph, necessarily negative. This is known as the {\it Delsarte bound}, as it was proved by Delsarte in \cite{Delsarte} for strongly regular graphs, and later  generalized by Godsil to distance-regular graphs (see \cite[Proposition 4.4.6]{BCN}). We will refer to a clique in $\Gamma$ as a {\em Delsarte clique} if its order is equal to $1+\frac{k}{-\theta_{\min}}$, in other words if achieves equality in the Delsarte bound.

The following lemma shows why we prefer generic parameters over standard ones in relation to Delsarte cliques.

\begin{lem} \label{delsarte stuff}
Let $\Gamma$ be a strongly regular graph with generic parameters $(s, t, \alpha)$, and suppose that $\Gamma$ contains a Delsarte clique $C$. Then $|C| = s + 1$, and $|\Gamma(x) \cap C| = \alpha$ for every vertex $x \in V(\Gamma) \backslash C$. In particular, $s$, $t$, and $\alpha$ are positive integers.
\end{lem}

{\bf Proof:} 
Note that $k = s(t+1)$, $\mu = (t+1)\alpha$, and $\theta_{\min} = -t-1$. The lemma then follows immediately from \cite[Lemma~8.6.6]{SRG}.
\qed




\subsection{Partial linear spaces}\label{sec:pls}

An {\em incidence structure} is a tuple $(\mathcal{P}, \mathcal{L}, \mathcal{I})$ where $\mathcal{P}$ and $\mathcal{L}$ are non-empty disjoint sets and $\mathcal{I} \subseteq \mathcal{P} \times \mathcal{L}$. The elements of $\mathcal{P}$ and $\mathcal{L}$ are called {\em points} and {\em lines}, respectively. If $(p, \ell) \in \mathcal{I}$, then we say that $p$ is {\em incident} with $\ell$, or that $p$ is on the line $\ell$. 

A \emph{partial linear space} is an incidence structure such that each pair of distinct points are both incident with at most one line. The \emph{point graph} or \emph{collinearity graph} $\Gamma$ of a partial linear space $(\mathcal{P}, \mathcal{L}, \mathcal{I})$ is the graph with vertex set $\mathcal{P}$, and two distinct points are adjacent if and only if they are on a common line.

General strongly regular graphs need not be related to partial linear spaces, but many are the point graphs of partial linear spaces. Knowledge that a graph is associated to such a structure gives access to powerful results on the structure of the graph, and therefore a sufficient condition for a strongly regular graph to be of this type is of great value. The most important result on this was proved by Metsch in \cite{Metsch}.  

\begin{thm}[{\rm\hspace{1sp}cf.\cite[Result 2.1]{metsch1999characterization}, \cite{Metsch}}]\label{metsch}
Let $\Gamma$ be an amply regular graph with parameters $(v,k, \lambda, \mu)$.
Assume that there exists a positive integer $\sigma$ such that the following two conditions are satisfied:
\begin{itemize}
\item[$(1)$] $(\sigma+1)(\lambda +1) - k > (\mu-1){\sigma+1 \choose 2}.$
\item[$(2)$] $\lambda+1  > (\mu -1)(2\sigma -1)$.
\end{itemize}
Define a line as a maximal clique with at least $\lambda+2 - (\mu -1)(\sigma-1)$ vertices.
Then $X= (V(\Gamma), {\mathcal L}, \in)$ is a partial linear space, where $\mathcal L$ is the set of all lines, $\Gamma$ is the point graph of $X$, and the symbol $\in$ means that the relation required for the incidence structure is given by inclusion. Moreover, each vertex lies on at most $\sigma$ lines.
\end{thm}

 Note that the condition $(1)$ of Theorem~\ref{metsch} implies that each vertex lies on at most $\sigma$ lines because it means that there are no $(\sigma+1)$-claws (see, for example \cite[Lemma 1.6]{kgpp}).

Motivated by Theorem~\ref{metsch}, we introduce the SPLS($\sigma$) property, first defined in \cite{lv2024bounding}. 

Let $\Gamma$ be a strongly regular graph with standard parameters $(v,k,\lambda,\mu)$, and let $\sigma \geq 2$ be an integer. We say that $\Gamma$ has the {\it SPLS($\sigma$)} property if the following conditions hold:  
\begin{itemize}
\item[(1)] $\Gamma$ is the point graph of a partial linear space $(V(\Gamma), \mathcal{L}, \in)$, where the set of lines $\mathcal{L}$ are the maximal cliques with at least $\lambda+2-(\mu-1)(\sigma-1)$ vertices.  
\item[(2)] $\lambda \geq (2\sigma-1)(\mu-1)$.  
\item[(3)] Each vertex lies on at most $\sigma$ lines.  
\end{itemize}

These conditions come directly from the conditions in Theorem \ref{metsch}. In particular, if $\Gamma$ satisfies the conditions of that theorem for some $\sigma$, then it has the SPLS($\sigma$) property. We note further that if $\Gamma$ has the SPLS($\sigma$) property, then a lower bound on the size of each line is 
\begin{equation} \label{splssigma}
\lambda+2-(\mu-1)(\sigma-1) \geq \sigma(\mu-1)+2.
\end{equation}

\subsection{Partial geometries}\label{sec:partialgeo}

A \emph{partial geometry} $pg(s,t,\alpha)$ is a partial linear space with the property that every line contains exactly $s+1$ points, every point lies on exactly $t+1$ lines, and given a line $\ell$ and a point $x \notin \ell$, there are exactly $\alpha$ lines containing $x$ and intersecting $\ell$. (In the literature, one also meets the notation $pg(K,R,T)$, where $K=s+1$, $R=t+1$, and $T=\alpha$.) We will assume that $s,t,\alpha>0$.  

Given a $pg(s,t,\alpha)$, denoted by $X$, we may dualize it; in other words, we may create a new partial geometry $Y$ whose points are the lines of $X$ and whose lines are the points of $X$, with a point lying on a line in $Y$ precisely when the corresponding line contains the corresponding point in $X$. In this process we obtain a $pg(t,s,\alpha)$.  

The following important result on partial geometries was proved by Neumaier.

\begin{thm} \cite[Theorem 4.5]{neumaier-m} \label{fakea}
If a $pg(s,t,\alpha)$ satisfies $\alpha \leq t-1$, then
$$
s \leq (t - \alpha + 1)^2 (2\alpha - 1),
$$
and equality implies $\alpha = 1$ or $t = 2\alpha$.
\end{thm}

The dualization process described earlier allows us to interchange $s$ and $t$ in this result, and in fact we have stated the dual of the result as it was originally stated in \cite{neumaier-m}.

The point graph of a partial geometry $pg(s,t,\alpha)$ is automatically strongly regular, and has standard parameters $(s + 1 + (s + 1) s t/\alpha, s (t + 1), s + t (\alpha - 1) - 1, (t + 1) \alpha)$. We will say that a strongly regular graph which is the point graph of a partial geometry is {\it geometric}. Note that cliques of order $s+1$ in the point graph of a partial geometry $pg(s,t,\alpha)$ are Delsarte cliques; therefore, each line in the partial geometry corresponds to a Delsarte clique in the point graph.

However, it is possible for a strongly regular graph to have standard parameters of the above form but which is not the point graph of a partial geometry; it is even possible that no partial geometry with the corresponding parameters exist (an example is the Cameron graph, see the introduction of \cite{guo2020nonexistence} for a discussion of this). 

Let $\Gamma$ be a strongly regular graph with standard parameters $(v, k, \lambda, \mu)$. We say that $\Gamma$ has \emph{geometric parameters} $(s, t, \alpha)$ if there exist positive integers $s, t, \alpha$ with $\alpha \leq s+1$ and $\alpha \leq t+1$ that satisfy the following conditions:

\begin{equation} \label{geopar}
\begin{gathered} 
\hspace{-0.2cm}k = s(t+1),\\
\hspace{1.2cm}\lambda = s + t(\alpha - 1 ) - 1, \\
\hspace{-0.25cm}\mu = (t + 1)\alpha. \\
\end{gathered}    
\end{equation}

\vspace{14pt}

The following lemma describes the relationship between the generic and geometric parameters, which follows directly from the definition.

\begin{lem}\label{genericgeo}
Let $\Gamma$ be a primitive strongly regular graph with generic parameters $(s, t, \alpha)$. Then $\Gamma$ admits geometric parameters if and only if $s$, $t$, and $\alpha$ are positive integers satisfying $\alpha \leq s+1$ and $\alpha \leq t+1$.
\end{lem}



\section{Strongly regular graphs that arise as the point graph of a partial linear space} \label{sec: partial lin space}

In this section, we apply the methods developed in \cite{guo2020nonexistence} to a general setting, in order to provide a collection of lemmas which are required for Theorem \ref{mainclassification}. Delsarte cliques, and the vertices contained in them, are of principal importance throughout. 

Let $\Gamma$ be a strongly regular graph with generic parameters $(s, t, \alpha)$. Recall that $\theta_{\min} = -t-1$. We will assume that $\Gamma$ is the
point graph of a partial linear space $X = (V(\Gamma), \mathcal{L}, \in)$, where the set of lines $\mathcal{L}$ consists
of a certain collection of maximal cliques of $\Gamma$.  For each $x \in V(\Gamma)$, let $\tau(x)$ denote the number of lines in the partial linear space $X$ containing $x$. The following lemma, giving us a lower bound on $\tau(x)$, is important. 

\begin{lem} \label{lowboundtau}
$\tau(x) \geq t+1$ for all $x$.
\end{lem}

{\bf Proof:} Recall from Lemma \ref{delsarte stuff} that the maximal size of a line is $s+1$, and by the definition of a partial linear space any two distinct lines containing $x$ can only intersect at $x$. Therefore, any collection of $t'+1<t+1$ lines containing $x$ would contain less than $s(t+1)$ vertices in $\Gamma(x)$, however there are $k = s(t+1)$ vertices in $\Gamma(x)$, so there must be some vertices in $\Gamma(x)$ not contained in any of these $t'+1$ lines. Since every edge is contained in a unique line, we see that there must be a line containing $x$ which is not among our collection of $t'+1$ lines. It follows that there must be at least $t+1$ lines containing $x$. \qed

\vspace{14pt} 

It follows also from this argument that, if $\tau(x) = t+1$ for some $x$, then every line containing $x$ must be a Delsarte clique, containing exactly $s + 1$ vertices. We will call such a vertex (with $\tau(x) = t+1$) a {\it Delsarte vertex}, and let $V_D$ denote the set of Delsarte vertices in $V(\Gamma)$. Note that if all vertices are Delsarte vertices, then $\Gamma$ is geometric.


We also will need the following estimate on the number of lines.

\begin{lem} \label{numline}
We have $\min(g,|V_D|) \geq |V(\Gamma)| - |{\cal L}|$, where $g$ is the multiplicity of the smallest eigenvalue $-t-1$ of $\Gamma$.
\end{lem}

{\bf Proof:} We recall that ${\cal L}$ denotes the set of lines in the partial linear space. Let $M$ denote the {\it line-vertex incidence matrix} of $\Gamma$. That is, $M$ is a $|{\cal L}| \times |V(\Gamma)|$ matrix whose rows are indexed by the elements of ${\cal L}$ and columns are indexed by the vertices of $\Gamma$, with $M_{ij} = 1$ when vertex $x_j$ is contained in line $L_i$ and 0 otherwise. Since every edge is contained in a unique line, we will have

$$
M^T M = A + diag(\tau(x_1),\tau(x_2), \ldots, \tau(x_n)).
$$
where $V(\Gamma)=\{x_1, \ldots, x_n\}$ and $diag(\cdot)$ denotes the diagonal matrix with the given values on the diagonal. It is clear that  $rank(M) \leq \min(|{\cal L}|, |V(\Gamma)|)$.


We know that the smallest eigenvalue of $A$ is $-t-1$ of multiplicity $g$, and the eigenvalues of $diag(\tau(x_1),\tau(x_2), \ldots, \tau(x_n))$ are simply the quantities $\{\tau(x_j)\}$, all of which are at least $t+1$ by Lemma \ref{lowboundtau}. 

Suppose $v$ is a unit 0-eigenvector of $M^TM$. Then $v^T(A+diag(\tau(x_1),\tau(x_2), \ldots, \tau(x_n)))v = 0$. However, we must also have $v^TAv \geq -t-1$ and $v^Tdiag(\tau(x_1),\tau(x_2), \ldots, \tau(x_n))v \geq t+1$. These inequalities must therefore be equalitites, and thus $v$ must necessarily be a $(t+1)$-eigenvector of $diag(\tau(x_1),\tau(x_2), \ldots, \tau(x_n))$ and a $(-t-1)$-eigenvector of $A$. We see that $|V_D|$, which is the dimension of the $(t+1)$-eigenspace of $diag(\tau(x_1),\tau(x_2), \ldots, \tau(x_n))$, is at least as big as the dimension of the 0-eigenspace of $M^TM$, which is $|V(\Gamma)|-rank(M) \geq |V(\Gamma)|-|{\cal L}|$. The same is true of $g$, the dimension of the $(-t-1)$-eigenspace of $A$. The result follows. \qed

\vspace{14pt}

The following lemma combines the previous lemma with a simple upper bound on $|{\cal L}|$ in order to give a lower bound on $|V_D|$ in the case $\mu \geq 2$. In fact, the estimate does hold when $\mu=1$, but gives no information since $\sigma \geq 2$. As such, we need a new argument that holds when $\mu = 1$, and this is given by the subsequent Lemma \ref{mu1} below.

\begin{lem} \label{numDel}
Let $\Gamma$ be a strongly regular graph with standard parameters $(v,k,\lambda,\mu)$. If $\mu \geq 2$ and $\Gamma$ has the SPLS($\sigma$) property, then 

$$
|V_D| \geq \Big(1-\frac{\sigma}{\sigma(\mu - 1) + 2}\Big) |V(\Gamma)|.
$$

In particular, if $\mu \geq 3$ then $|V_D| \geq \frac{1}{2} |V(\Gamma)|$, and if $\mu=2$ then $|V_D| > 2$.
\end{lem}

{\bf Proof:} Note that $|{\cal L}|$ can be upper bounded by summing the number of lines containing each vertex and dividing by the minimum size of a line. In other words 

$$|{\cal L}| \leq \frac{\sum_{x \in V(\Gamma)}\tau(x)}{(\lambda+2)-(\mu-1)(\sigma-1)}  \leq \frac{|V(\Gamma)| \sigma}{\sigma(\mu-1)+2},
$$
where the second inequality follows from $(\ref{splssigma})$ and the fact that each vertex lies on at most $\sigma$ lines.
By Lemma~\ref{numline}, we know that $|V_D|\geq|V(\Gamma)|-|\cal L|$, and hence we have
$$|V_D| \geq \Big(1-\frac{\sigma}{\sigma(\mu - 1) + 2}\Big) |V(\Gamma)|.$$

For the final statement, if $\mu \geq 3$ then the fraction $\frac{\sigma}{\sigma(\mu - 1) + 2}$ is less than $\frac{1}{2}$, so the result follows. If $\mu=2$ then $1-\frac{\sigma}{\sigma(\mu - 1) + 2} = \frac{2}{\sigma+2}$. The conditions $\sigma \geq 2$ and $\lambda \geq (2\sigma - 1)(\mu-1) = 2\sigma - 1$, required for the SPLS($\sigma$) property, imply that $|V(\Gamma)| > \lambda+2 > \sigma+2$. Thus, $|V_D| \geq \frac{2}{\sigma+2}|V(\Gamma)| > 2$. \qed

\vspace{12pt}

The fact that this lemma does not apply to the case $\mu = 1$ is not a serious obstacle, since it turns out that when $\mu = 1$ the graph in question can never be geometric, and there is a simple bound  on $\lambda$ in terms of $k$. In particular, we have the following result.

\begin{lem}\label{mu1}
Let $\Gamma$ be a strongly regular graph with standard parameters $(v,k,\lambda,\mu)$. If $\mu = 1$, then $\Gamma$ is not geometric and $(\lambda + 1)(\lambda + 2) \leq k$.
\end{lem}

{\bf Proof.}
Suppose that $\Gamma$ is a geometric strongly regular graph, and fix a vertex $x \in V(\Gamma)$. By Lemma \ref{delsarte stuff}, every vertex not adjacent to $x$ has at least one neighbour on each of the $\frac{k}{\lambda+1}$ lines containing $x$, and these neighbours are all distinct. Therefore, $\mu \geq \frac{k}{\lambda+1} \geq 2$. This yields our first conclusion. 

For the second, assume $\mu=1$ and $(\lambda + 1)(\lambda + 2) > k$. This implies that $\Gamma$ has the
SPLS($\lambda + 1$) property, since in this case we clearly have $\lambda+1 \geq 2$, and all the conditions of Theorem \ref{metsch} are satisfied for $\sigma= \lambda+1$. Now we have

$$
|{\cal L}| \leq \frac{\sum_{x \in V(\Gamma)} \tau(x)}{(\lambda+2) - (\mu-1)(\sigma-1)} \leq \frac{|V(\Gamma)|(\lambda+1)}{(\lambda+2)} <|V(\Gamma)|.
$$

Hence, there exists at least one Delsarte vertex in $\Gamma$ (by Lemma~\ref{numline}).
Since $\mu = 1$, for each vertex $x$ of $\Gamma$, any two nonadjacent vertices in $\Gamma(x)$ cannot have a common neighbor in $\Gamma(x)$. Also, clearly each vertex in $\Gamma(x)$ has exactly $\lambda$ neighbors in $\Gamma(x)$ as $\Gamma(x)$ is $\lambda$-regular. Hence, $\Gamma(x)$ is a disjoint union of $\frac{k}{\lambda + 1}$
cliques, each of order $\lambda + 1$. In other words, each vertex is contained in
$\frac{k}{\lambda + 1}$ maximal cliques of order $\lambda + 2$. Therefore, if a Delsarte
vertex exists, every vertex in $\Gamma$ is a Delsarte vertex, and thus $\Gamma$ is
geometric. This contradicts our earlier conclusion, and the result follows. \qed

The following lemma imposes strong restrictions on the generic parameters.

\begin{lem}\label{alphaint}
Let $\Gamma$ be a strongly regular graph with generic parameters $(s, t, \alpha)$. 
If $\Gamma$ contains two distinct Delsarte vertices, then it admits geometric parameters. In particular, $\alpha \leq t+1$.
\end{lem}

\noindent
{\bf Proof.}
By Lemma~\ref{delsarte stuff}, $s$, $t$, and $\alpha$ are positive integers. Furthermore, it follows from (\ref{classparsrg}) that $\alpha \leq s$. 

Let $x$ and $y$ be two distinct Delsarte vertices in $\Gamma$. If $d(x, y) = 1$, then there exists a Delsarte clique $C$ that contains $y$ but not $x$. Since $C$ is a maximal clique, there must be a vertex $z \in C$ that is not adjacent to $x$. If $d(x, y) = 2$, let $z = y$. In either case, we can find a vertex $z$ with $d(x, z) = 2$ such that there exists a Delsarte clique $C$ containing $z$ and $x$ has at least one neighbor in $C$.

By Lemma~\ref{delsarte stuff}, $x$ has exactly $\alpha$ neighbors in $C$. 
Since $x$ is a Delsarte vertex, there are $ t + 1$ lines through $x$. Let $C_1, C_2, \ldots, C_{t+1}$ denote these lines, each of which is a Delsarte clique.
Note that $x$ has exactly $\alpha$ neighbors in $C$. Furthermore, each clique $C_i$ intersects $C$ in at most one vertex for $i = 1, 2, \ldots, t+1$. This implies that $\alpha \leq t + 1$.  Consequently, by Lemma~\ref{genericgeo}, $\Gamma$ admits geometric parameters. \qed

\section{Proof of Theorem \ref{mainclassification}} \label{sec: proof}

In this section, we will finally prove our main result. In Theorem~\ref{splsprop}, we establish a sufficient condition for a strongly regular graph with generic parameters $(s,t,\alpha)$ to have the \SPLS property. Lemma~\ref{pseudogeo} shows that such graphs possess geometric parameters.  
In Theorem \ref{geometric} we will show, under weak conditions, that a strongly regular graph with geometric parameters is in fact geometric. At the end of the section, we will combine everything to prove Theorem \ref{mainclassification}.

\begin{thm}
\label{splsprop}

Let $\Gamma$ be a strongly regular graph with generic parameters $(s,t,\alpha)$, where $t \geq 1$. Suppose that $\Gamma$ has standard parameter $\mu \geq 2$ and that 

\begin{align*}
s > \max \Big\{&\frac{8}{3}(t+1)(t\alpha+\alpha-1) - \frac{2}{3}(t\alpha+\alpha-1) - 4t(\alpha-1),\\
&\frac{8}{3}(t+1)(t\alpha+\alpha-1) - \frac{5}{3}(t\alpha+\alpha-1) - t(\alpha-1)\Big\}.
\end{align*}

Then $\Gamma$ has the SPLS($\lceil \frac{4t+1}{3} \rceil$) property.

\end{thm}

\vspace{12pt}

\textbf{Proof:} By definition of the SPLS($\sigma$) property, we only need to show that the conditions (1) and (2) in Theorem \ref{metsch} are all satisfied for $\sigma = \lceil \frac{4t+1}{3} \rceil$. Clearly, $\sigma=\lceil\frac{4t+1}{3}\rceil\geq 2$ as $t\geq1$. From $\lambda+1 = s+t(\alpha-1)$, $\mu-1 = t\alpha+\alpha-1$, and $k = s(t+1)$, the conditions in Theorem~\ref{metsch} can be written as

\begin{align}
 s &> \frac{\sigma + 1}{2(\sigma - t)}((t\alpha+\alpha-1)\sigma - 2t(\alpha-1)), \label{eq:a} \\
 s &> (t\alpha+\alpha-1)(2\sigma - 1) - t(\alpha-1). \label{eq:b}
\end{align}


Let $\sigma_0 := \frac{4t+1}{3}$ be a real number. Then $\sigma_0 \leq \sigma \leq \sigma_0 + \frac{2}{3}$, and hence $\frac{\sigma+1}{\sigma-t} \leq \frac{\sigma_0+1}{\sigma_0-t} = 4$.

As $(t\alpha+\alpha-1)\sigma - 2t(\alpha-1) > 0$ and $t\alpha+\alpha-1 = \mu-1\geq1$ always hold, to show that conditions \eqref{eq:a} and \eqref{eq:b} hold, it suffices to show the following two inequalities hold:
\begin{align*} 
    &s > 2(t\alpha+\alpha-1)(\sigma_0+\frac{2}{3})-4t(\alpha-1) = \frac{8}{3}(t+1)(t\alpha+\alpha-1) - \frac{2}{3}(t\alpha+\alpha-1) - 4t(\alpha-1), \\
     &s > (t\alpha+\alpha-1)(2\sigma_0 + \frac{4}{3}-1)-t(\alpha-1) = \frac{8}{3}(t+1)(t\alpha+\alpha-1) - \frac{5}{3}(t\alpha+\alpha-1) - t(\alpha-1).
\end{align*}
Therefore, if $\begin{aligned}[t]
s > \max \Big\{&\frac{8}{3}(t+1)(t\alpha+\alpha-1) - \frac{2}{3}(t\alpha+\alpha-1) - 4t(\alpha-1),\\
&\frac{8}{3}(t+1)(t\alpha+\alpha-1) - \frac{5}{3}(t\alpha+\alpha-1) - t(\alpha-1)\Big\},
\end{aligned}$

\noindent then $\Gamma$ has the SPLS($\lceil \frac{4t+1}{3} \rceil$) property by Theorem \ref{metsch}. \qed

\vspace{12pt}

\begin{lem}
\label{pseudogeo}
Let $\Gamma$ be a strongly regular graph with generic parameters $(s,t,\alpha)$, where $t \geq 1$. Suppose that $\Gamma$ has standard parameter $\mu \geq 2$. 
If $\Gamma$ has the $SPLS(\sigma)$ property for some $\sigma$, then it admits geometric parameters. In particular, $\alpha \leq t+1$.

\end{lem}

\textbf{Proof:} Lemma~\ref{numDel} implies that there are at least two Delsarte vertices, and the result then follows from Lemma~\ref{alphaint}.
\qed

\begin{thm}
\label{geometric}
    Let $\Gamma$ be a strongly regular graph with generic parameters $(s,t,\alpha)$, where $t \geq 2$. Suppose that $\Gamma$ has standard parameter $\mu \geq 2$ and has the $SPLS(\sigma)$ property for some $\sigma \leq 2t$. If $s \geq \frac{5}{2}t(\bal)$, then $\Gamma$ is a geometric strongly regular graph.
\end{thm}

\textbf{Proof: }Assume $\Gamma$ satisfies the conditions of the theorem, but is not geometric. By Lemma~\ref{pseudogeo}, $t$ and $\alpha$ are positive integers with $\alpha \leq t+1$. Moreover, $\mu = (t+1)\alpha \geq t+1 \geq 3$. 

We will let $V^c_D$ denote the set of vertices which are not Delsarte vertices; that is, $V^c_D = V(\Gamma)-V_D = \{x\in V(\Gamma)~|~\tau(x)>t+1\}$. As $\Gamma$ is not geometric, $V^c_D \neq \phi$. For any vertex $x$ of $\Gamma$, let $\tau_D(x)$ denote the number of Delsarte cliques containing $x$; note that $\tau_D(x) \leq \tau(x)$, with equality precisely when $x$ is a Delsarte vertex. Fix integer $\gamma$ with $0 \leq \gamma \leq  t-1$, and set $\delta =  t+1-\gamma$, note that $2 \leq \delta \leq t+1$. We now separate the analysis into two cases, the first when there is $x \in V^c_D$ with $\tau_D(x) \geq \gamma+1$, and the second when there is no such $x$. We will obtain a bound on $s$ in both cases, and later will show how to choose $\gamma$ in order to optimize this bound.

{\bf Case 1:} Suppose there is $x \in V^c_D$ with $\tau_D(x) \geq \gamma+1$. Then we will be able to show that
\begin{equation*}
    s \leq \binom{\delta}{2}(\bal)-(\delta-2)\delta(\alpha-1)\text{, for $2 \leq \delta \leq t+1$}.
\end{equation*}

To prove this, let $\tau(x) = t+2+\kappa$ for some $\kappa \geq 0$. Since $\Gamma$ has the \SPLS property, we must have $\tau(x) \leq \sigma \leq 2t$. Note that $\tau_D(x) \geq \gamma+1$ implies that there are at least $\gamma+1$ Delsarte cliques containing $x$. We fix $\gamma+1$ Delsarte cliques containing $x$, and we count the vertices of $\Gamma(x)$ which lie outside of the $\gamma+1$ Delsarte cliques containing $x$. 
Let $C_1, C_2, \ldots, C_{\gamma+1}, D_1, D_2, \ldots, D_{\tau(x)-\gamma-1}$ be the lines containing $x$, where $C_1, C_2, \ldots, C_{\gamma+1}$ are Delsarte cliques. Note that $|C_i| = s+1$ for $1 \leq i \leq \gamma+1$.
Let $D_i' = D_i \setminus \{x\}$, and choose a vertex $y_i \in D_i'$ for $1 \leq i \leq \tau(x)-\gamma-1$.
Let $T = D_1' \cup D_2' \cup \cdots \cup D_{\tau(x)-\gamma-1}'$, and define 
$T_i = (\Gamma(y_i) \cup \{y_i\}) \cap T$ for $1 \leq i \leq \tau(x)-\gamma-1$.
By the inclusion-exclusion principle, we have

\begin{align*}
    k-s(\gamma+1) &= |T| = \left|\bigcup_{i=1}^{\tau(x)-\gamma-1} T_i\right|
\geq \sum_{i=1}^{\tau(x)-\gamma-1} |T_i|
\;-\;
\sum_{1 \leq i < j \leq \tau(x)-\gamma-1} |T_i \cap T_j| \\
&\geq (\tau(x)-\gamma-1)(\lambda+1-(\gamma+1)(\alpha-1))-\binom{\tau(x)-\gamma-1}{2}(\mu-1).
\end{align*}

As $k=s(t+1), \lambda=s+t(\alpha-1)-1, \mu=(t+1)\alpha, \tau(x)-\gamma-1=\delta+\kappa$ and $\gamma+1=t+2-\delta$, we have
\begin{align} \label{stars}
    s(t-\gamma) &\geq (\delta+\kappa)(s + (\delta - 2)(\alpha - 1)) - \binom{\delta+\kappa}{2}(\bal). 
\end{align}


We wish to replace the right-hand side of \eqref{stars} by $\delta(s + (\delta - 2)(\alpha - 1))-\binom{\delta}{2}(\bal)$, and this is justified by the following sequence of equivalent inequalities. Note that clearly it is true when $\kappa=0$, so we may assume that $\kappa>0$.

\begin{align*}
(\delta+\kappa)(s& + (\delta - 2)(\alpha - 1)) - \binom{\delta+\kappa}{2}(\bal) \geq \delta(s + (\delta - 2)(\alpha - 1))-\binom{\delta}{2}(\bal) \nonumber \\
&\iff \kappa(s + (\delta - 2)(\alpha - 1)) \geq \Big(\binom{\delta+\kappa}{2}-\binom{\delta}{2}\Big)(\bal) \nonumber\\
&\iff \kappa(s + (\delta - 2)(\alpha - 1)) \geq \frac{\kappa^2+2\delta \kappa-\kappa}{2}(\bal) \nonumber\\
&\iff s + (\delta - 2)(\alpha - 1) \geq \frac{\kappa+2\delta-1}{2}(\bal) \nonumber\\
&\iff s \geq (2-\delta)(\alpha-1) + \frac{\kappa+2\delta -1}{2}(\bal).
\end{align*}

As $\tau(x) \leq 2t$, $\kappa=\tau(x)-t-2 \leq t-2$ and $s \geq \frac{5}{2}t(\bal)$, the above holds.\\
In other word, we have 

$$s(t-\gamma-\delta) \geq \delta(\delta-2)(\alpha - 1)-\binom{\delta}{2}(\bal).$$

Recall that $\delta = t+1-\gamma$, and substituting in this final expression for $\delta$ on the left yields $-s \geq (\delta-2) \delta (\alpha - 1)-\binom{\delta}{2}(\bal)$. Multiplying by $-1$ gives our bound in Case 1.

\vspace{12pt}

{\bf Case 2:} Suppose $\tau_D(x) \leq \gamma$ for all $x \in V^c_D$. Then we will show that \\
\begin{equation} \label{daggers}
    s \leq \frac{(\bal)t(\alpha-1) \delta + (t+1)(t+\alpha)+1}{(\bal)\delta(\delta-1)-(t+1)t}\text{, for $2 \leq \delta \leq t+1$}. 
\end{equation}

To see this, note that by Lemma~\ref{numDel}, $|V^c_D| \leq \frac{\sigma}{\sigma(\mu-1)+2}|V(\Gamma)|<\frac{1}{\mu-1}|V(\Gamma)|$. We will now establish a lower bound on $|V^c_D|$. Consider $\Gamma(x)$ for any $x \in V^c_D$. If $y\in \Gamma(x) \setminus V^c_D$, then $y$ is a Delsarte vertex. By our assumption, there are at most $s\gamma$ such vertices in $\Gamma(x)$. Therefore \\
\begin{equation*}
    |\Gamma(x) \cap V^c_D| \geq k-s\gamma = s(t+1-\gamma)=s\delta.
\end{equation*}

Next, we want to count the number of vertices $z \in \Gamma_2(x)\cap V^c_D$. By the same arguement as before, each $y \in \Gamma(x) \cap V^c_D$ must also have at least $s \delta$ neighbours in $V^c_D$. Out of these, we must remove $x$ and those already in $\Gamma(x)$. There are at most $\lambda + 1 = s + t(\alpha-1)$ of those.

Finally, for each $z \in \Gamma_2(x) \cap V^c_D$, it is counted at most $\mu$ times via elements of $\Gamma(x)$. Thus, a lower bound on $|V^c_D|$ is given by
\begin{equation*}
    |V^c_D| \geq 1+s \delta + \frac{s\delta(s\delta-s-t(\alpha-1))}{(t+1)\alpha} > \frac{s\delta(s\delta-s-t(\alpha-1))}{(t+1)\alpha}.
\end{equation*}
Note that $|V(\Gamma)|=1+k+\frac{k(k-\lambda-1)}{\mu}=1+s(t+1)+\frac{s t(s-\alpha+1)}{\alpha}$.

By $|V^c_D| < \frac{|V(\Gamma)|}{\mu-1}=\frac{|V(\Gamma)|}{t\alpha+\alpha-1}$, we have
\begin{align*}
    & \frac{s \delta (s \delta - s - t(\alpha-1))}{(t+1)\alpha} < \frac{1}{\bal} \left( \frac{s t(s - \alpha + 1)}{\alpha} + s(t+1)+1 \right) \\
    \iff\hspace{0.2cm} &s\delta (s \delta - s - t(\alpha-1))(\bal) \\
    & < (t+1)ts(s-\alpha + 1)+(t+1)^2\alpha s + (t+1)\alpha \\
    \overset{\text{divide by }s}{\iff} &(\bal)(\delta-1)\delta s - (\bal)t(\alpha-1)\delta \\
    & < (t+1)ts - (t+1)t(\alpha-1) + (t+1)^2 \alpha  +\frac{(t+1)\alpha}{s} \\
    \iff\hspace{0.2cm} & ((\bal)(\delta-1)\delta - (t+1)t)s \\
    & < (\bal)t(\alpha-1)\delta + (t+1)(\alpha + t)+\frac{(t+1)\alpha}{s} \\
    & < (\bal)t(\alpha-1)\delta + (t+1)(\alpha + t)+1.
\end{align*}

The final inequality is true because $s \geq \frac{5}{2}t(\bal)$ and $t \geq 2, \alpha \geq 1$. Dividing through to isolate $s$ on the left gives the bound \eqref{daggers} we have claimed for Case 2. 

\vspace{12pt}

It is evident that exactly one of these two cases must hold, so combining the estimates we have

\begin{align}
    \label{eqn2}
    \nonumber s &\leq \max \mathopen{\Biggl\{} 
    \binom{\delta}{2}(\bal)-(\delta-2)\delta(\alpha - 1), \\
    &\phantom{\leq \max \mathopen{\Biggl\{}} 
    \frac{(\bal)t(\alpha-1)\delta + (t+1)(t+\alpha)+1}{(\bal)\delta(\delta-1)-(t+1)t} 
    \Biggr\}.
\end{align}

Recall that $\gamma$ was chosen arbitrarily between $0$ and $t-1$, and so this bound holds for any $\delta$ with $2 \leq \delta \leq t+1$. We will now show that $\delta$ can be chosen so that the right hand side of \eqref{eqn2} is less than $\frac{5}{2}t(t\alpha+\alpha-1)$.

If $\alpha=1$ and $t=2$, it is easy  to check that the right side of \eqref{eqn2} is 2, which less than  $\frac{5}{2}t(t\alpha+\alpha-1)=10$, for $\delta=2$. Therefore, we may assume that either ($\alpha > 1$ and $t\geq2$) or ($\alpha = 1$ and  $t\geq3$).

Take $\delta=\lfloor \sqrt{5t} \rfloor$. It is easy to check that $2 \leq \delta \leq t+1$ holds, since $t \geq 2$. This implies that $\binom{\delta}{2}(\bal) - (\delta-2)\delta (\alpha - 1) <\frac{\delta^2}{2}(\bal) \leq \frac{5}{2}t(\bal)$, so we need only check the second expression on the right of \eqref{eqn2}.

Note that $\sqrt{5t}-1 < \delta \leq \sqrt{5t}$. We have

\begin{align*}
    &(\bal)\delta(\delta-1)-(t+1)t \\
    &> (\bal)(\sqrt{5t}-1)(\sqrt{5t}-2)-(t+1)t \\
    &= (\bal)(5t+2-3\sqrt{5t})-(t+1)t \\
    &= (\bal)+\Bigl((\bal)(5t+1-3\sqrt{5t})-(t+1)t\Bigl).
\end{align*}
We will show $(\bal)(5t+1-3\sqrt{5t})-(t+1)t \geq 0$. To show this, note that  if $\alpha>1$ and $t\geq2$, then
\begin{align*}
    (\bal)(5t+1-3\sqrt{5t})-(t+1)t \geq (2t+1)(5t+1-3\sqrt{5t})-(t+1)t \\
    =9t^2+6t+1-6t\sqrt{5t}-3\sqrt{5t} >0.
\end{align*}
If $\alpha=1$ and $t \geq 3$, it is easy to check $t(5t+1-3\sqrt{5t})-(t+1)t \geq 0$. Hence, the second term of \eqref{eqn2} has an upper bound
\begin{align*}
    \frac{(\bal)t(\alpha-1)\delta + (t+1)(t+\alpha)+1}{\bal} &\leq t(\alpha-1)\delta + \frac{(t+1)(t+\alpha)+1}{\bal} \\
    &\leq t(\alpha-1)\sqrt{5t}+\frac{(t+1)(t+\alpha)+1}{\bal}.
\end{align*}
As $\mu = (t+1)\alpha \geq 3$, we have $\bal \geq 2$, and therefore

\begin{align*}
    t(\alpha-1)\sqrt{5t}+\frac{(t+1)(t+\alpha)+1}{\bal} 
    &\leq t(\alpha-1)\sqrt{5t}+\frac{(t+1)(t+\alpha)+1}{2} \\
    &< \frac{5}{2}t(\bal) \quad \text{holds}.
\end{align*}

We see that, in either case, taking $\delta = \lfloor \sqrt{5t} \rfloor$ shows that the right side of \eqref{eqn2} is less than $\frac{5}{2}t(\bal)$. This contradicts our assumption, which was that $s \geq \frac{5}{2}t(\bal)$. This completes the proof, for we have shown that, under the given conditions, all vertices in $\Gamma$ are Delsarte vertices. In other words, $\Gamma$ is a geometric strongly regular graph. \qed

\vspace{12pt}

\if2 
The following can be considered a corollary of the previous result, but since it is our main result we list it as a theorem.

\begin{thm}
    \label{final result 1}
    Let $\Gamma$ be a strongly regular graph with $\mu \geq 2$ and smallest eigenvalue $-m$, where $m \geq 3$ is an integer. If $\mu \neq m(m-1), m^2$ then

    \begin{equation}
        \label{eqn3}
        \lambda \leq \frac{8}{3}m (\mu-1)-\frac{2}{3}\mu+3m-\frac{10}{3}.
    \end{equation}
\end{thm}
\fi 

We now are in a position to prove our main result.

\vspace{12pt}

\textbf{Proof of Theorem \ref{mainclassification}: } If $m=1$, then by \cite[Corollary 3.5.4]{BCN}, $\Gamma$ is the disjoint union of complete graphs. If $m = 2$, then $\Gamma$ is known to belong to a short list of graphs (see \cite[p. 5]{SRG})), and the result is easily checked for those graphs. We may therefore assume $m \geq 3$. 

Suppose $\mu = 1$. Then we only need to show that $(iii)$ holds.
Assume $\lambda > f(m, \mu) = 3m - 4$, i.e., $\lambda \geq 3m - 3$.
By \cite[Lemma 3.2]{-m-koolenbang}, we have $k < m(\lambda + m)$.
Then $(\lambda + 1)(\lambda + 2) > m(\lambda + m) > k$, which contradicts Lemma \ref{mu1}. This proves $(iii)$ when $\mu = 1$, so we may assume $\mu \geq 2$ in the following.

Let $\Gamma$ has generic parameters $(s, t, \alpha)$. Then $t=m-1\geq2$. As $\lambda = s+t(\alpha-1)-1$, $\mu=(t+1)\alpha$ and $m=t+1$, $\lambda \leq f(m,\mu)$ is equivalent to $s\leq\frac{8}{3}(t+1)(\bal)-\frac{2}{3}t \alpha-\frac{2}{3}\alpha +3t +\frac{2}{3}-t(\alpha-1)$.

Assume $s>\frac{8}{3}(t+1)(\bal)-\frac{2}{3}t \alpha-\frac{2}{3}\alpha +3t +\frac{2}{3}-t(\alpha-1)$. Note that $\alpha = \frac{\mu}{t+1}>0$. \\
It is easy to check that 

\begin{align*}
    s > {\rm max} \biggl\{ &\frac{8}{3}(t+1)(\bal)-\frac{2}{3}(\bal)-4t(\alpha-1), \\
    &\frac{8}{3}(t+1)(\bal)-\frac{5}{3}(\bal)-t(\alpha-1) \biggl\}.
\end{align*}

Therefore, by Theorem \ref{splsprop}, $\Gamma$ has the SPLS($\lceil \frac{4t+1}{3} \rceil$) property. It is easy to check that $s \geq \frac{5}{2}t(\bal)$ and $\lceil \frac{4t+1}{3} \rceil \leq 2t$. Therefore, by Theorem \ref{geometric}, $\Gamma$ is a geometric strongly regular graph. Now, $(V(\Gamma), \mathcal{L}, \in)$ is a partial geometry $pg(s, t, \alpha)$, where $\mathcal{L}$ are all Delsarte cliques of $\Gamma$. If $\mu=m(m-1)$(resp. $\mu = m^2$), then we have $\alpha = t$ (resp. $\alpha = t+1$), and hence we are in case $(i)$ (resp. case $(ii)$) of the theorem. If neither of these hold, then $\alpha \leq t-1$, and Theorem \ref{fakea} therefore implies

\begin{equation*}
    s \leq(t-\alpha+1)^2(2\alpha-1)<2t^2(\alpha-1)<\frac{5}{2}t(\bal). 
\end{equation*}

This is a contradiction, and proves the theorem. \qed

\vspace{12pt}

\begin{rem} \label{finalthoughts}

\begin{enumerate}

\item From the proof of Theorem~\ref{mainclassification}, we know that if $\mu = 1$, then $\lambda
\leq 3m - 4$. In fact, by applying a similar argument, we can derive the following result: Let
$\Gamma$ be a strongly regular graph with standard parameter $\mu = 1$ and smallest eigenvalue
$-m$. Then, we have $\lambda < 1.62m$.

\item \cite[Theorem 2.3]{Metsch} also showed that if a strongly regular graph has geometric parameters, then it is geometric, under slightly stronger constraints, by different methods from
those we employed in Theorem \ref{geometric}. We remark that the first two conditions in \cite[Theorem 2.3]{Metsch} are exactly the same as the two conditions in Theorem \ref{metsch}.

\item By employing more refined calculations, \cite[Lemma 2.4]{Metsch}  provides a slightly
improved bound for constructing partial linear spaces in comparison to our Theorem \ref{splsprop}, under the
assumption that $\alpha \geq \frac{3}{11}(t+1)+\frac{21}{11}$. Using the same computational techniques, we
can also improve our bound.

\item In \cite{guo2020nonexistence}, it was shown that for $\alpha = 1$, one has $s \leq t \lfloor \frac{8}{3}t + 1
\rfloor$. By our Theorem \ref{splsprop}, Theorem \ref{geometric} and Theorem \ref{fakea}, we have $s \leq \frac{8}{3}t^2
+ 2t$. The result in \cite{guo2020nonexistence} is slightly better because there the authors carried out a more detailed discussion by considering different cases, which we have not done here.

\item Comparing our main result, Theorem \ref{mainclassification}, with the result of Neumaier,
Theorem \ref{srgneumaier}, we note that the leading term in our bound is $ \frac{8}{3} m \mu $, while the
leading term in his bound is $ \frac{1}{2} m^2 \mu $. In general, our bound is better; however,
due to a larger constant term in our expression, it only becomes better when $ m $ is sufficiently
large. Specifically, our bound is better than his when $ m \geq 6 $. Apart from a few
exceptional cases, his bound is almost always better when $ m \leq 5 $.

\end{enumerate}
\end{rem}

\section{Comparing our bound with the bound of Babai and Wilmes}\label{compBWBound}

In this section, we compare our bound on the parameter $\lambda$ with the result obtained by Babai and Wilmes for primitive strongly regular graphs.

In \cite{BabaiWilmes2016}, Babai and Wilmes proved that the parameter $\lambda$ satisfies  
\begin{equation}\label{BWBound}
    \lambda + 1 < \max\left\{4\sqrt{2v}, \frac{6}{\sqrt{13} - 1}\sqrt{k(\mu - 1)}\right\}.
\end{equation}

We will show that the bound given in Theorem~\ref{mainclassification},  
\begin{equation}\label{ourbound}
    \lambda \leq \frac{8}{3} m (\mu - 1) - \frac{2}{3} \mu + 3m - \frac{10}{3},
\end{equation}
is stronger than the one stated above. To prove this, it suffices to show that inequality~\eqref{ourbound} implies inequality~\eqref{BWBound}. Since we are interested in the case where the smallest eigenvalue $-m$ is an integer, and noting that when $m = 1$ the graph is a disjoint union of complete graphs, we may therefore assume that $m \geq 2$.

Since $\frac{6}{\sqrt{13} - 1} \approx 2.303$, the inequality $\lambda + 1 < \frac{6}{\sqrt{13} - 1} \sqrt{k(\mu - 1)}$ is approximately equivalent to $(\lambda + 1)^2 < 2.303^2 k(\mu - 1)$.  We proceed by considering three cases.

(1) If $k > \frac{4}{5} m(\lambda+1)$ and $\lambda + 1 < 4m(\mu - 1)$, then $(\lambda + 1)^2 < 4m(\mu - 1) (\lambda + 1)  < 5k(\mu - 1) < 2.303^2 k(\mu - 1)$.

(2) If $k > \frac{4}{5} m (\lambda+1)$ and $\lambda + 1 \geq 4m(\mu - 1)$, then by inequality~\eqref{ourbound}, we obtain $4m(\mu - 1) \leq \frac{8}{3} m (\mu - 1) - \frac{2}{3} \mu + 3m - \frac{7}{3} < 3m\mu$. Therefore, $\mu \leq 3$. Since $k > \frac{4}{5} m (\lambda+1)$, we have $k - \lambda - 1 \geq \frac{3}{5}  (\lambda+1)$, and hence $v = 1 + k + \frac{k(k - \lambda - 1)}{\mu} \geq 1 + k + \frac{3k(\lambda + 1)}{5\mu} > \frac{4}{25} m (\lambda + 1)^2 \geq \frac{8}{25} (\lambda + 1)^2$. Thus, we obtain $\lambda + 1 < 4 \sqrt{2v}$.

(3) Assume $k \leq \frac{4}{5} m (\lambda+1)$. By equation~\eqref{eigeq}, we have $k = \mu + m \lambda + m^2 - \mu m$, and hence $\frac{1}{5} m \lambda \leq \frac{4}{5} m + \mu m - m^2 - \mu$. Therefore, $\lambda + 1 < 5(\mu + 1 - m) \leq 5 (\mu - 1)$. It follows that $(\lambda + 1)^2 < 5 (\mu - 1)(\lambda + 1) \leq 5k(\mu - 1) < 2.303^2 k(\mu - 1)$.


\appendix
\section{Appendix: Applications to orthogonal arrays} \label{sec: orth array}

Our methods provide another approach to proving Bruck's Completion Theorem, as we now explain. An {\em orthogonal array} $\mathrm{OA}(m,n)$ is an $m \times n^2$ array
with entries from an $n$-element set, such that all $n^2$ ordered pairs formed
by any two rows are distinct. An $\mathrm{OA}(m,n)$ is equivalent to a set of
$m - 2$ mutually orthogonal Latin squares.

Given an $\mathrm{OA}(m,n)$, we define the {\em Latin square graph}
$\mathrm{LS}_m(n)$ as follows: the vertices are the $n^2$ columns of the
orthogonal array, and two vertices are adjacent if they agree in one coordinate
position. This graph is strongly regular with standard parameters $(v, k, \lambda, \mu) =
(n^2,\, m(n-1),\, (m-1)(m-2) + n - 2,\, m(m-1))$, with second largest eigenvalue
$n - m$ and smallest eigenvalue $-m$. A necessary condition for the existence of an orthogonal array
$\mathrm{OA}(m,n)$ is that $m \leq n + 1$. An $\mathrm{OA}(n+1,n)$ is
called a {\em full orthogonal array}. The {\em deficiency} of an
$\mathrm{OA}(m,n)$ is defined as $\delta = n - m + 1$, representing the
number of rows missing from a full orthogonal array. The complement of
$\mathrm{LS}_m(n)$ has the same parameters as $\mathrm{LS}_\delta(n)$.

We remark that certain references use the term {\it Latin square graph} for any strongly regular graph with standard parameters of the form $(n^2,\, m(n-1),\, (m-1)(m-2) + n - 2,\, m(m-1))$, regardless of whether it corresponds to an orthogonal array.
However, we will stick to the definition given above.

The result of Bruck, Theorem \ref{srgneumaier} $(i)$ (cf.\ \cite{bruck1963finite}), implies that if $n > f(\delta)$,
where $f$ is a fixed polynomial of degree 4, then the complement is itself a
Latin square graph. In this case, the two graphs can be combined to obtain a full
orthogonal array. Metsch \cite[Corollary 2.6]{Metsch} later improved
Bruck’s bound by providing a polynomial $f$ of degree 3.

Following the proof strategy of Bruck and Metsch, and applying Theorems \ref{splsprop}
and \ref{geometric}, we obtain a bound similar to that of Metsch, namely the following
corollary.

\begin{cor} \label{OA}
Let $m$ and $n$ be positive integers, and define $\delta = n - m + 1$. Suppose there
exists an orthogonal array $\mathrm{OA}(m,n)$. If $n > \frac{8}{3} \delta^3 - \frac{16}{3} \delta^2 + 2\delta + \frac{2}{3}$,
then the orthogonal array can be extended to a full orthogonal array.    
\end{cor}

Here we provide an outline of the proof; for full details, one may refer to \cite{bruck1963finite} or \cite[Corollary 2.5]{Metsch}. \\

{\bf Proof (sketch):} If $\delta = 1$, the statement is clearly true \cite[Corollary 2.5]{Metsch}, so we
may assume $\delta \geq 2$. Assume we have an orthogonal array
$\mathrm{OA}(m,n)$, and let $\delta = n - m + 1$. Then we can construct the Latin
square graph $\mathrm{LS}_m(n)$. The complement of $\mathrm{LS}_m(n)$, which we denote
by $\Gamma$, has the same parameters as the Latin square graph
$\mathrm{LS}_\delta(n)$, which has smallest eigenvalue $-\delta$, and where $n >
\frac{8}{3} \delta^3 - \frac{16}{3} \delta^2 + 2\delta + \frac{2}{3}$. If $\delta = 2$, by
the classification of strongly regular graphs with smallest eigenvalue $-2$ (cf. \cite[p. 5]{SRG}), $\Gamma$ is $\mathrm{LS}_2(n)$, which is geometric. If $\delta \geq 3$, then
by Theorems \ref{splsprop} and \ref{geometric}, we can also show that the graph $\Gamma$ is geometric. The
set of all lines of $\Gamma$ can be partitioned into $\delta$ classes, each forming a
partition of $V(\Gamma)$. For each such partition, we assign the same label to all
vertices lying on the same line. In this way, we naturally obtain a new line of the
orthogonal array. Since there are $\delta$ such partitions, we obtain $\delta$ new lines,
thereby extending the orthogonal array to a full orthogonal array.

\begin{rem}
\begin{enumerate}
\item In \cite{bruck1963finite}, Bruck shows that every orthogonal array $\mathrm{OA}(m,n)$ with $m > n -
O(n^{1/4})$ can be extended to a full orthogonal array.

\item Metsch \cite[Corollary 2.6]{Metsch} and our Corollary \ref{OA} show that every orthogonal array
$\mathrm{OA}(m,n)$ with $m > n - O(n^{1/3})$ can be extended to a full orthogonal array.

\item It seems likely the current methods would not allow one to obtain the following conclusion, which was mentioned by Bruck in \cite{bruck1963finite}:
every orthogonal array $\mathrm{OA}(m,n)$ with $m > n - O(n^{1/2})$ can be extended to
a full orthogonal array.

\item Essentially, our Corollary \ref{OA} and Metsch's result \cite[Corollary~2.6]{Metsch} are the same.
However, Metsch uses a more refined calculation when constructing the partial linear space (see
\cite[Lemma 2.4]{Metsch}), compared to our Theorem \ref{splsprop}. As a result, his bound is slightly better. If we adopted
the same technique, we would obtain the same bound as he did.

\item In \cite{Metsch}, Metsch also applied Bruck’s Completion Theorem
to other topics in finite geometry.
\end{enumerate}

\end{rem}

\section*{Acknowledgements}

J.H. Koolen is partially supported by the National Natural Science Foundation of China (No. 12471335) and the Anhui Initiative in Quantum Information Technologies (No. AHY150000). J. Park was supported by the National Research Foundation of Korea (NRF) grant funded by the Korea government (MSIT) (RS-2024-00356153). G. Markowsky would like to thank Harrison Choi for his assistance, and J.H. Koolen would like to thank Sanming Zhou and Melbourne University for a generous invitation which led to the results in this paper. We thank the anonymous referee for pointing out the connection between our results and the graph isomorphism problem.

\bibliographystyle{plain}
\bibliography{conf}
\end{document}